\documentclass{amsproc}

\usepackage[swedish,english]{babel}
\usepackage[utf8]{inputenc}
\usepackage[T1]{fontenc}
\usepackage{amsmath,amssymb,amsfonts}
\usepackage{url}
\usepackage[dvipsnames,usenames]{color}
\usepackage{array}
\usepackage[pdftex]{graphicx}
\usepackage{accents}

\newcommand{\R}{\mathbb R}
\newcommand{\N}{\mathbb N}

\newcommand{\C}{\mathbb C}

\newcommand{\cont}{\mathcal C}

\renewcommand{\Im}{\mathrm{Im}}
\renewcommand{\Re}{\mathrm{Re}}
\def\E{{\mathrm{e}}}
\def\di{\partial} 
\def\til{\widetilde}
\def\I{\mathfrak{i}}
\newcommand{\diff}{\mathrm{d}}
\renewcommand{\bar}{\overline}

\newcommand{\Dom}{\mathrm{Dom}}
\newcommand{\dd}{\mathcal D}
\newcommand{\sd}{\mathcal S}
\newcommand{\ed}{\mathcal E}
\newcommand{\supp}{\mathrm{supp}}
\newcommand{\Leb}{\mathrm{L}}
\newcommand{\Loc}{\Leb^1_{\mathrm{loc}}}

\newcommand{\Ran}{\mathrm{Ran}}
\newcommand{\holo}{holomorphic }

\newcommand{\ie}{\textit{i.e.}\/ } % in other words, that is to say % no comma after in British English
\newcommand{\eg}{\textit{e.g.}\/ } % for example % no comma after in British English
\newcommand{\cf}{\textit{cf.}\/ } % compare with

\renewcommand{\vec}[1]{\accentset{\rightharpoonup}{#1}} % vector symbol
\allowdisplaybreaks

\theoremstyle{definition} 
\newtheorem{define}{Definition}[section]
\newtheorem{example}[define]{Example}
\newtheorem{remark}[define]{Remark}

\theoremstyle{plain}
\newtheorem{lemma}[define]{Lemma}
\newtheorem{thm}[define]{Theorem}
\newtheorem{prop}[define]{Proposition}

\numberwithin{equation}{section}

\begin{document}

\title[pseudo-passive causal operators of slow growth]{On pseudo-passive causal operators of slow growth}

\author{Mitja Nedic}
\address{Mitja Nedic, Department of Mathematics, Stockholm University, SE-106 91 Stockholm, Sweden, orc-id: 0000-0001-7867-5874}
\curraddr{}
\email{mitja@math.su.se}
\thanks{The author is supported by the Swedish Foundation for Strategic Research, grant nr. AM13-0011.}


\begin{abstract}
In this paper, we study a class of convolution operators on the space of distributions that enlarge the well-studied class of passive operators. In this larger class, we are able to associate, to each operator, a \holo function in the right half-plane with a specific constraint on its range, determined by the operator. Afterwards, we investigate whether the properties of causality and slow growth hold automatically in our larger class of convolution operators. Finally, an alternative class of convolution operators is also considered.
\end{abstract}

\maketitle

\section{Introduction}

Convolution operators have applications in many areas of science, \eg mechanics, thermodynamics and electromagnetism. In electromagnetic theory, many processes can be represented as a relation between the input and output signal of some process or device called a \emph{system}. As the signals entering and leaving a system are mathematically often modeled by distributions, we may consider the system to be represented by an operator on the space of distributions. Under some basic assumptions such as linearity and continuity, \cf Remark \ref{rem:schwartz_characterization}, any operator on the space of distributions may be written as a convolution with a fixed distribution \cite{Schwartz1952,Zemanian1968a}. Standard examples of such systems and their modeling convolution operators can be found among electrical circuits, \eg \cite[Sec. 5.1]{BernlandEtal2011} and \cite[pg. 315]{Vladimirov1979}.

The theory of passive systems, or \emph{passive operators}, is the theory of such convolution operators, which also satisfy the condition of passivity, \cf Definition \ref{def:passive_operator}. In electromagnetics, this condition is interpreted as the system being unable to produce its own energy. The usefulness of this theory is dependent on a classical result which states that every convolution operator satisfying the condition of passivity may be represented in terms of a \holo function
on the right half-plane having non-negative real part, \cf Theorem \ref{thm:passive_operators}.

Outside of electromagnetics, the classical theory of passive operators has been considered in the scalar-valued setting, \eg \cite{Nussenzveig1972,WohlersBeltrami1965,Zemanian1965}, in the matrix-valued setting, \eg \cite{Vladimirov1979,YoulaEtal1959,Zemanian1968a}, and in the operator-valued setting, \eg \cite{Zemanian1972}. More recently, greater focus has been placed on the study of infinite-dimensional systems, \eg \cite{GuiverLogemannOpmeer2017,Staffans2002}, and applications of operator-valued setting, \eg \cite{HanygaSeredynska2008}.

The physics interpretation of the theory of passive operators is also its drawback. Hence, one would like to be able to adapt the classical theory to a class of operators which satisfy a more general condition than passivity, one which could be interpreted as the electromagnetic system either producing some of its own energy, or being supplied energy form an outside source. We are thus interested in identifying a more general condition than the classical condition of passivity that would allow us generalize the main result of the classical theory. 

In this paper, we define a condition called \emph{pseudo-passivity} and successfully translate some of the most important results of the classical theory to this larger setting. More precisely, we define a class of convolution operators we call \emph{pseudo-passive causal operators of slow-growth}, \cf Definition \ref{def:non_passive_operator}, and associate, to each such operator, a \holo function in the right half-plane, \cf Theorem \ref{thm:range_of_W}. Furthermore, we investigate the relationship between the conditions of pseudo-passivity, causality and slow growth, \cf Theorem \ref{thm:implications}.

The structure of the paper is as follows. In Section \ref{sec:background}, we briefly recall the main results of the theory of passive operators and introduce the class of pseudo-passive operators. Section \ref{sec:Laplace} is then devoted to establishing the correspondence between these operators and certain \holo functions in the right half-plane, while Section \ref{sec:causality_slow_growth} investigates whether the conditions of causality and slow growth are automatically satisfied within the class of pseudo-passive operators. Finally, in Section \ref{sec:scattering_operators}, we consider an analogous generalization of the theory of scattering passive operators, \cf Definitions \ref{def:scattering_passive_operator} and \ref{def:scattering_non_passive_operator}, while some concluding statements are presented in Section \ref{sec:summary}.

\section{Background}\label{sec:background}

Let us first set the notation that we will use throughout the paper and recall the necessary bits and pieces of distribution theory \cite{Zemanian1965}. The right-half plane is denoted by $\C_+ := \{z \in \C~|~\Re[z] \geq 0 \}$ and we introduce the abbreviations $\Leb^1 := \Leb^1(\R,\C)$ and $\Loc := \mathrm{L}^1_{\mathrm{loc}}(\R,\C)$. Let also $\dd := \cont_0^\infty(\R,\C)$ and $\dd'$ be the usual spaces of test functions and distributions, $\sd$ and $\sd'$ the spaces of Schwartz functions and Schwartz distributions, and $\ed$ and $\ed'$ the spaces of smooth functions and compactly supported distributions, respectively. We define the inclusion $\iota\colon \Loc \hookrightarrow \dd'$ as $\iota\colon f \mapsto T_f$, where the distribution $T_f$ is defined, for any test function $\varphi \in \dd$, as
$$\langle T_f,\varphi \rangle := \int_\R f(t)\varphi(t)\diff t.$$

We say that a distribution $U$ is \emph{non-zero} on an open subset $\Omega \subseteq \R$ if there exists a test function $\varphi \in \dd$ with $\supp(\varphi) \subseteq \Omega$ such that $\langle U,\varphi \rangle \neq 0$. Furthermore, a distribution $U$ is called \emph{non-negative definite} if
\begin{equation*}\label{eq:non_neg_definite}
\left\langle U, t \mapsto \int_{-\infty}^\infty\bar{\varphi(\tau)}\varphi(\tau - t)\diff\tau\right\rangle \geq 0
\end{equation*}
for any test function $\varphi \in \dd'$ \cite[pg. 304]{Zemanian1965}.

The \emph{support} of a distribution $U$ is defined as
$$\supp(U) := \{t \in \R~|~U \text{ is non-zero in every open neighbourhood of } t\}.$$
An important subspace of distributions is the space of \emph{right-sided distributions} $\dd'_\mathrm{r}$, which consists of all distributions $U \in \dd$ such that there exists a number $t \in \R$ (depending on $U$) for which $\supp(U) \subseteq [t,\infty)$.

Another important object in the theory of operators on distributions is the \emph{Laplace transform} \cite[Sec. 8]{Zemanian1965}. For a function $f$ with $\supp(f) \subseteq [t,\infty)$ for some $t \in \R$ and $\tau \mapsto \E^{-c\:\tau}f(\tau) \in \Leb^1$ for some $c \in \R$, we define its Laplace transform $\mathcal{L}(f)$, for $s \in \C$ with $\Re[s] \geq c$, as
$$\mathcal{L}(f)(s) := \int_t^\infty f(\tau)\E^{-s\:\tau}\diff\tau.$$
If $U$ is a distribution, assuming that $U \in \dd'_\mathrm{r}$ and that $(\tau \mapsto \E^{-c\:\tau})U \in \sd'$ for some $c \in \R$, we define its Laplace transform $\mathcal{L}(U)$, for $s \in \C$ with $\Re[s] \geq c$, as
$$\mathcal{L}(U)(s) := \langle (\tau \mapsto \E^{-c\:\tau})U, \xi \mapsto \rho(\xi)\E^{-(s-c)\xi} \rangle,$$
where $\rho$ is any smooth function with support bounded on the left such that $\rho \equiv 1$ on some neighbourhood of $\supp(U)$.

A fundamental concept in the theory of distributions is the \emph{convolution}. First, we note that for two distributions $U,V \in \dd'$ there exists precisely one distribution $V \otimes U \in (\cont^\infty_0(\R^2,\C))'$ with the property that
$$\langle V \otimes U, (t_1,t_2) \mapsto \varphi_1(t_1)\varphi_2(t_2) \rangle = \langle V, \varphi_1 \rangle\langle U,\varphi_2 \rangle.$$
The distribution $V \otimes U$ is called the \emph{direct product} of $V$ and $U$ and can be defined in two equivalent ways, either as
$$\langle V \otimes U, \psi \rangle :=  \langle V, t_1 \mapsto \langle U, t_2 \mapsto \psi(t_1,t_2)\rangle\rangle$$
or as
$$\langle V \otimes U, \psi \rangle := \langle U, t_2 \mapsto \langle V, t_1 \mapsto \psi(t_1,t_2)\rangle\rangle,$$
where $\psi \in \cont^\infty_0(\R^2,\C)$ is any test function (in two variables), \cf \cite{Nussenzveig1972} and \cite[Sec. 5.2]{Zemanian1965}. We may now define the convolution $V*U$ as
$$\langle V * U, \varphi \rangle := \langle V \otimes U, (t_1,t_2) \mapsto \varphi(t_1+t_2) \rangle,$$
where $\varphi \in \dd$ is any test function. We note that the direct product of two distributions always exists, but the convolution of two distributions does not necessarily exist. This follows from the observation that the function $(t_1,t_2) \mapsto \varphi(t_1+t_2)$ never has compact support (unless $\varphi \equiv 0$), yielding restrictions on the supports of $U$ and $V$. However, it is not a problem to define the convolution of two distributions if both are in $\sd'$, or if both are in $\dd'_\mathrm{r}$, or as long as at least one of them has compact support.

Consider now a convolution operator $R$ on $\dd'$ with defining distribution $Y \in \dd'$, \ie $R = Y*$, where $\Dom(R)$ is taken as the largest possible domain of definition \cite{Zemanian1965,Zemanian1968a}. By the definition of the convolution between distributions, we may be certain that $\Dom(R)$ always contains at least all compactly supported distributions. Furthermore, a convolution operator $R$ is called \emph{real} if it maps real distributions to real distributions, where, we recall, a distribution $U \in \dd'$ is called real if $\langle U,\varphi \rangle \in \R$ for any test function $\varphi \in \dd$ taking only real values.

Let us now recall the standard definition of passive operators \cite{Zemanian1965}.

\begin{define}\label{def:passive_operator}
A convolution operator $R = Y *$ is called a \emph{passive operator} if it satisfies the following three conditions:
\begin{itemize}
    \item[(a)]{\emph{passivity} on $\iota(\dd)$: for any number $t \in \R$ and any test function $\varphi \in \dd$ it holds that $R(\iota(\varphi)) \in \iota(\Loc)$ and
\begin{equation}\label{eq:passivity}
\Re\left[\int_{-\infty}^t\bar{\varphi(\tau)}\psi(\tau)\diff\tau\right] \geq 0,
\end{equation}
where $\psi := \iota^{-1}(R(\iota(\varphi)))$,}
    \item[(b)]{\emph{causality}: $\supp(Y) \subseteq [0,\infty)$,}
    \item[(c)]{\emph{slow growth}: $Y \in \sd'$.}
\end{itemize}
\end{define}

\begin{remark}
The appearance of the conditions of causality and slow growth in Definition \ref{def:passive_operator} turns out to be superfluous, \cf Section \ref{sec:causality_slow_growth}. Any convolution operator that satisfies the condition of passivity will automatically be causal and of slow growth. For these reasons, operators satisfying Definition \ref{def:passive_operator} are referred to as \emph{passive operators} instead of \emph{passive causal operators of slow growth}. When other types of convolution operators are also being considered, \cf Section \ref{sec:scattering_operators}, we refer to operators satisfying Definition \ref{def:passive_operator} as \emph{admittance passive operators}.
\end{remark}

\begin{remark}\label{rem:schwartz_characterization}
A convolution operator on $\dd'$ may be characterized via four conditions: single-valuedness, linearity, strong continuity and time-translation invariance, \cf \cite{Schwartz1959} and \cite[Thm. 2]{Zemanian1968a}.
\end{remark}

We note that condition of passivity requires that the operators $R$ maps distributions arising from test functions to distributions arising from locally integrable functions.

\begin{example}
A first example, let us check the operator $R = \delta_0'*$ is a passive operator.

If $\varphi \in \dd$ is any test function, we note first that
$$\delta_0' * T_\varphi = T_{\varphi'}$$
and, as such, it holds that
$$\iota^{-1}(\delta_0' * T_\varphi) = \varphi'.$$
Using integration by parts, we now calculate, for any $t \in \R$, that
$$\int_{-\infty}^t\bar{\varphi(\tau)}\varphi'(\tau)\diff \tau = |\varphi(t)|^2 - \int_{-\infty}^t\varphi(\tau)\bar{\varphi'(\tau)}\diff\tau.$$
Reorganizing the above equality yields
$$\Re\left[\int_{-\infty}^t\bar{\varphi(\tau)}\varphi'(\tau)\diff\tau\right] = \frac{1}{2}|\varphi(t)|^2 \geq 0,$$
showing that the operator $R$ is indeed passive.\hfill$\lozenge$
\end{example}

It is a well-known result that any passive operator can be uniquely described in terms of the Laplace transform $W := \mathcal{L}(Y)$ of its defining distribution $Y$. The conditions of causality and slow growth assure the existence of $W$ as a \holo function in the right half-plane, while the condition of passivity restricts the range of the function $W$. In particular, the following class of functions is to be considered \cite{Zemanian1965,Zemanian1968a}.

\begin{define}\label{def:PR_function}
A \holo function $p\colon \C_+ \to \C$ for which $\Re[p(s)] \geq 0$ for all $s \in \C_+$ and $p(s) \in \R$ for $s \in (0,\infty)$ is called a \emph{positive-real function}.
\end{define}

The following theorem now describes the correspondence between real passive operators and positive-real functions, \cf \cite[Thm. 10.4-1]{Zemanian1965} and \cite[Thm. 10.6-1]{Zemanian1965}.

\begin{thm}\label{thm:passive_operators}
Let $R = Y*$ be a real passive operator. Then, the Laplace transform of its defining distribution $Y$ exists and is a positive-real function. Conversely, for any positive-real function $W$, the operator $R := \mathcal{L}^{-1}(W)*$ is a real passive operator.
\end{thm}

\begin{remark}
Without the assumption that the operator $R= Y*$ is real, it still holds that the Laplace transform of its defining distribution $Y$ exists and is \holo function in the right half-plane having non-negative real part, \cf \cite[Thm. 10.4-1]{Zemanian1965}. Conversely, for any \holo function in the right half-plane having non-negative real part, its inverse Laplace transform $\mathcal{L}^{-1}(W)$ exists and can be used to define a convolution operator, \cf \cite[Thm. 10.6-1]{Zemanian1965}.
\end{remark}

The generalization of Definition \ref{def:passive_operator} that we are interested in is the class of convolution operators given by the following definition.

\begin{define}\label{def:non_passive_operator}
A convolution operator $R = Y *$ is called a \emph{pseudo-passive causal operator of slow growth} if it satisfies the following three conditions:
\begin{itemize}
\item[(a')]{\emph{pseudo-passivity} on $\iota(\dd)$: there exist a number $N \in \N_0$ and vectors $\vec{c},\vec{d} \in \R^{N+1}$, such that for any number $t \in \R$ and any test function $\varphi \in \dd$, it holds that $R(\iota(\varphi)) \in \iota(\dd)$ and
\begin{equation}\label{eq:admittance_non_passivity}
\Re\left[\int_{-\infty}^t\bar{\varphi(\tau)}\psi(\tau)\diff\tau\right] \geq \sum_{j = 0}^{N}\int_{-\infty}^t\left(c_j|\varphi^{(j)}(\tau)|^2 + d_j|\psi^{(j)}(\tau)|^2\right)\diff\tau,
\end{equation}
where $\psi := \iota^{-1}(R(\iota(\varphi)))$ (if $\vec{d} = \vec{0}$, the requirement that $R(\iota(\varphi)) \in \iota(\dd)$ may be weakened to $R(\iota(\varphi)) \in \iota(\Loc)$),}
\item[(b)]{\emph{causality}: $\supp(Y) \subseteq [0,\infty)$,}
\item[(c)]{\emph{slow growth}: $Y \in \sd'$.}
\end{itemize}
\end{define}

\begin{example}
The first example of a pseudo-passive causal operator of slow growth is $-\delta_0*$, \ie a convolution with the negative of the Dirac distribution. Indeed, since
$$\iota^{-1}(\delta_0 * T_\varphi) = \varphi$$
for any test function $\varphi \in \dd$, it holds, for any $t \in \R$, that
$$\Re\left[\int_{-\infty}^t\bar{\varphi(\tau)}\psi(\tau)\diff\tau\right] = \Re\left[\int_{-\infty}^t\bar{\varphi(\tau)}\iota^{-1}(-\delta_0 * T_\varphi)(\tau)\diff\tau\right] = -\int_{-\infty}^t|\varphi(\tau)|^2\diff\tau.$$
In other words, the operator $-\delta_0*$ satisfies the condition of pseudo-passivity for $N = 0$, $c_0 = -1$ and $d_0 = 0$.\hfill$\lozenge$
\end{example}

\begin{example}
Assume that an operator $R$ satisfies the condition of pseudo-passivity with equality for some number $N \in \N_0$ and vectors $\vec{c},\vec{d} \in \R^{N+1}$. Then, the operator $-R$ satisfies condition \eqref{eq:admittance_non_passivity} with equality for the same number $N$ and vectors $-\vec{c},-\vec{d} \in \R^{N+1}$.\hfill$\lozenge$
\end{example}

In the next two sections, we will now focus on the following questions. First, how does the condition of pseudo-passivity restrict the range of the Laplace transform of the operator, \cf Section \ref{sec:Laplace}, and second, are the conditions of causality and slow growth superfluous as in the passive case, \cf Section \ref{sec:causality_slow_growth}.

\section{The Laplace transform}\label{sec:Laplace}

It is clear that the Laplace transform of the defining distribution of a pseudo-passive causal operator of slow growth exists, as this is assured by conditions of causality and slow growth. Therefore, we may investigate how the condition of pseudo-passivity restricts the range of the Laplace transform.

First, let us adopt an existing lemma to suit the condition of pseudo-passivity, \cf \cite[Lem. 6]{Zemanian1963} and \cite[Lem. 1, pg. 305]{Zemanian1965}.

\begin{lemma}\label{lem:addmitance_non_passive_schwartz}
Let $R$ be a pseudo-passive causal operator of slow growth. Then, the condition of pseudo-passivity is also satisfied for all Schwartz functions.
\end{lemma}

\proof
By the condition of slow growth, it holds that $\iota(\sd) \subseteq \Dom(R)$. Take, therefore, any Schwartz function $\varphi \in \sd$ and let $\{f_n\}_{n \in \N} \subseteq \dd$ is a sequence of test functions that converges to $\varphi$ in the topology of $\sd$. Since the operator $R$ satisfies the condition of pseudo-passivity, we have $R(\iota(f_n)) \in \iota(\dd)$ for all $n \in \N$ and, by \cite[Thm. 5.7-1]{Zemanian1965}, it holds that $R(\iota(\varphi)) \in \iota(\sd)$. Introducing $g_n := \iota^{-1}(R(\iota(f_n)))$, we may conclude, by \cite[Lem. 1, pg. 305]{Zemanian1965}, that
$$\left|\Re\left[\int_{-\infty}^t\bar{\varphi(\tau)}\psi(\tau)\diff\tau\right]-\Re\left[\int_{-\infty}^t\bar{f_n(\tau)}g_n(\tau)\diff\tau\right]\right| \to 0$$
as $n \to \infty$.

On the other hand, the right-hand side of inequality \eqref{eq:admittance_non_passivity} for the functions $f_n$ clearly converge to the right-hand side of inequality \eqref{eq:admittance_non_passivity} for the function $\varphi$. Furthermore, since inequality \eqref{eq:admittance_non_passivity} holds for any function $f_n$, it must also hold for the function $\varphi$. This finishes the proof.
\endproof

The established process of deriving the restriction on the range of the Laplace transform of the defining distribution of a passive operator, \cf \cite[pp. 306--307]{Zemanian1965}, may now be adapted to suit our generalization.

Take, therefore, $\varphi \in \sd$. Then, by Lemma \ref{lem:addmitance_non_passive_schwartz} and its proof, we have that $R(\iota(\varphi)) \in \iota(\sd)$ and it holds for $\psi := \iota^{-1}(R(\iota(\varphi)))$ that
$$\psi(\xi) = \langle Y, t' \mapsto \varphi(\xi - t')\rangle.$$
Furthermore, it also holds for any $t \in \R$ that
\begin{multline*}
\Re\left[\int_{-\infty}^t \bar{\varphi(\tau)}\psi(\tau)\diff\tau\right] = \Re\left[\int_{-\infty}^t \bar{\varphi(\tau)}\langle Y, t' \mapsto \varphi(\tau - t')\rangle\diff\tau\right] \\
\geq \sum_{j = 0}^{N}\int_{-\infty}^t\left(c_j|\varphi^{(j)}(\tau)|^2 + d_j|\psi^{(j)}(\tau)|^2\right)\diff\tau.
\end{multline*}

Let, therefore, $t \in \R$ be arbitrary and take any $s \in \C_+ := \{z \in \C~|~\Re[z] > 0\}$. Choose also a function $\varphi \in \sd$ with the property that $\varphi(\xi) = \E^{s\xi}$ for all $\xi \in (-\infty,a)$, where $a > t$. Since $\supp(Y) \subseteq [0,\infty)$, it holds that
\begin{multline*}
\Re\left[\int_{-\infty}^t \bar{\varphi(\tau)}\langle Y, t' \mapsto \varphi(\tau - t')\rangle\diff\tau\right] = \Re\left[\int_{-\infty}^t \E^{\bar{s}\tau}\langle Y, t' \mapsto \E^{s(\tau-t')}\rangle\diff\tau\right] \\
= \Re[\langle Y, t' \mapsto \E^{-s t'}\rangle]\int_{-\infty}^t \E^{2\Re[s]\tau}\diff\tau = \Re[W(s)]\int_{-\infty}^t \E^{2\Re[s]\tau}\diff\tau.
\end{multline*}

On the other hand, for such a function $\varphi$, it holds that
$$\psi(\xi) = \langle Y, t' \mapsto \E^{s(\xi-t')}\rangle = \E^{s\xi}\langle Y,t' \mapsto \E^{-st'}\rangle = \E^{s\xi}W(s)$$
for any $\xi \in (-\infty,t]$. For the derivatives of the functions $\varphi$ and $\psi$, it similarly holds that
$$\varphi^{(j)}(\xi) = s^j\E^{s\xi} \quad\text{and}\quad \psi^{(j)}(\xi) = s^j\E^{s\xi}W(s)$$
for any $j \in \N$ and any $\xi \in (-\infty,t]$. 

As such, we calculate that
\begin{multline*}
\sum_{j = 0}^{N}\int_{-\infty}^t\left(c_j|\varphi^{(j)}(\tau)|^2 + d_j|\psi^{(j)}(\tau)|^2\right)\diff\tau  \\
= \sum_{j = 0}^{N}|s|^{2j}\left(c_j+|W(s)|^2d_j\right)\int_{-\infty}^t\E^{2\Re[s]\tau}\diff\tau,
\end{multline*}
yielding
$$\Re[W(s)]\int_{-\infty}^t \E^{2\Re[s]\tau}\diff\tau \geq \sum_{j = 0}^{N}|s|^{2j}\left(c_j+|W(s)|^2d_j\right)\int_{-\infty}^t\E^{2\Re[s]\tau}\diff\tau.$$
Since the integral $\int_{-\infty}^t \E^{2\Re[s]\tau}\diff\tau$ is positive for any $s \in \C_+$, we conclude that
\begin{equation}\label{eq:range_of_W}
\Re[W(s)] \geq \sum_{j = 0}^{N}|s|^{2j}\left(c_j+|W(s)|^2d_j\right).
\end{equation}
Additionally, if the operator $R$ was real, this translates to the Laplace transform $W$ as the property that $W(s) \in \R$ if $s \in (0, \infty)$. Thus, we arrive at the following theorem.

\begin{thm}\label{thm:range_of_W}
Let $R = Y*$ be a pseudo-passive causal operator of slow-growth satisfying condition \eqref{eq:admittance_non_passivity} for some number $N \in \N_0$ and some vectors $\vec{c},\vec{d} \in \R^{N+1}$. Then, the Laplace transform $W := \mathcal{L}(Y)$ of its defining distribution exists and is a \holo function in the right half-plane satisfying condition \eqref{eq:range_of_W}. Furthermore, if the operator $R$ is assumed to be real, it holds that $W(s) \in \R$ if $s \in (0, \infty)$.
\end{thm}

\begin{remark}
It remains open whether all \holo functions in the right half-plane satisfying condition \eqref{eq:range_of_W} can be realized as the Laplace transform of the defining distribution of some pseudo-passive causal operator of slow growth. In the passive case, such a statement for the class of positive-real functions is made possible by the existence of an integral representation formula for which the inverse Laplace transform can be explicitly calculated, \cf \cite[Thm. 10.5-1]{Zemanian1965} and \cite[Thm. 10.6-1]{Zemanian1965}.
\end{remark}

\begin{example}\label{ex:cd_non_passive}
Let us consider the geometry of the case $N = 0$ of inequality \eqref{eq:range_of_W} in detail, \ie we are investigating the geometry of the range of a \holo function $W\colon \C_+ \to \C$ for which there exists numbers $c,d \in \R$, such that the inequality
\begin{equation}\label{eq:cd_non_passive_W}
\Re[W(s)] \geq c + d|W(s)|^2
\end{equation}
holds for any $s \in \C_+$. If this is the case, we have that $\Ran(W) \subseteq A$, where the set $A$ is defined as
$$A := \{\sigma \in \C~|~\Re[\sigma] \geq c+d|\sigma|^2\} = \{ (x,y) \in \R^2~|~x \geq c+dx^2+dy^2)\}.$$
Geometrically, the set $A$ is bounded by a circle or line in $\R^2$, with the precise picture being the following.

If $c = d = 0$, we have $A = \C_+ \cup \I\R$, as expected. Similarly, if $c \neq 0$, but $d = 0$, we get the half-plane $\{\sigma \in \C~|~\Re[\sigma] \geq c\}$. Therefore, it remains to consider the case $d \neq 0$ and $c \in \R$, where the inequality $x \geq c+dx^2+dy^2$ may be rewritten as
$$-\tfrac{c}{d} \geq (x^2 - \tfrac{x}{d}) + y^2 \quad\text{or}\quad \tfrac{1-4cd}{4d^2} \geq (x - \tfrac{1}{2d})^2 + y^2$$
if $d > 0$ and as
$$-\tfrac{c}{d} \leq (x^2 - \tfrac{x}{d}) + y^2 \quad\text{or}\quad \tfrac{1-4cd}{4d^2} \leq (x - \tfrac{1}{2d})^2 + y^2$$
if $d < 0$. Thus, the shape of the set $A$ is characterized in terms of the parameters $c$ and $d$ in the following way.
\begin{itemize}
\item[(i)]{We have $A = \C$ if $d < 0$ and $1-4cd \leq 0$.}
\item[(ii)]{We have $A$ being the area outside and including the circle $\tfrac{1-4cd}{4d^2} = (x - \tfrac{1}{2d})^2 + y^2$ if $d < 0$ and $1-4cd > 0$.}
\item[(iii)]{We have $A$ equal to the half-plane $\{(x,y) \in \R~|~x \geq c\}$ if $d = 0$.}
\item[(iv)]{We have $A$ being the area inside and including the circle $\tfrac{1-4cd}{4d^2} = (x - \tfrac{1}{2d})^2 + y^2$ if $d > 0$ and $1-4cd > 0$.}
\item[(v)]{We have $A = \{\tfrac{1}{2d}\} \subseteq \C$ if $d > 0$ and $1-4cd = 0$.}
\item[(vi)]{We have $A$ being empty if $d > 0$ and $1-4cd < 0$.}
\end{itemize}

In Figure \ref{fig:range_of_W} below, we see visualizations of the set A for different parameters $c$ and $d$. The parameters $(c,d) = (-1,-1)$ (top left) fall into case (i), the parameters  $(c,d) = (0,-1)$ (top centre), $(c,d) = (2,-1)$ (top right) and $(c,d) = (-\frac{1}{8},-\frac{1}{10})$ (middle left) all fall into case (ii), the parameters $(c,d) = (0,0)$ (middle centre) fall into case (iii), the parameters $(c,d) = (-2,\frac{1}{8})$ (middle right), $(c,d) = (0,\frac{1}{3})$ (bottom left) and $(c,d) = (0,1)$ (bottom centre) all fall into case (iv), and the parameters $(c,d) = (1,1)$ (bottom right) fall into case (vi).\hfill$\lozenge$
\end{example}

\begin{figure}[!ht]
\includegraphics[width=0.8\linewidth]{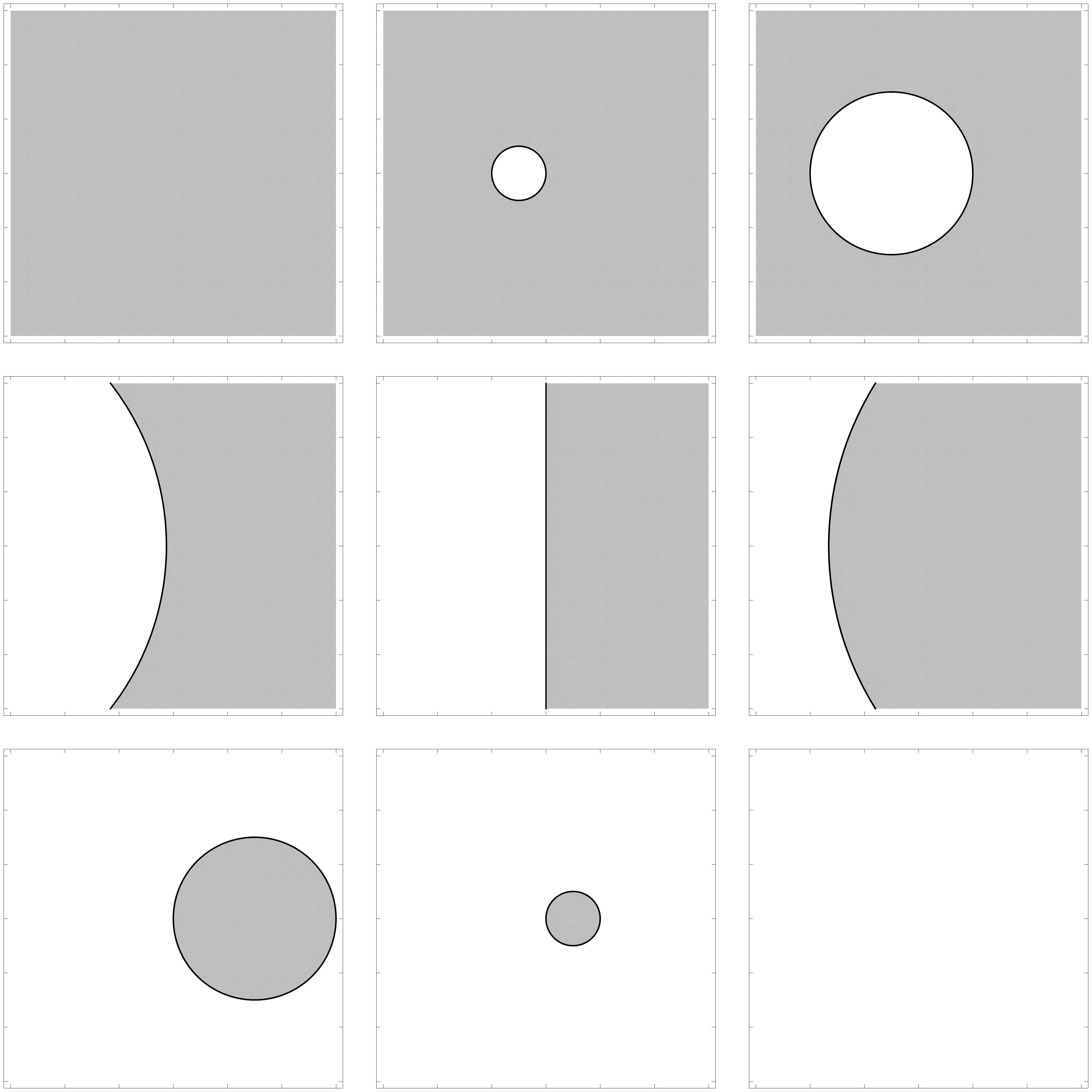}
\caption{Visualizations of the set $A$. The plot area is always the square $[-3,3]^2 \subseteq \R^2$.}
\label{fig:range_of_W}
\end{figure}

\section{Causality and slow growth}\label{sec:causality_slow_growth}

As already mentioned in Section \ref{sec:background}, any convolution operator that satisfies the condition of passivity \eqref{eq:passivity} on $\iota(\dd)$ must be given by a Schwartz distribution with support only on the positive half-line.

As we will soon see, this generalizes at least to some convolution operators satisfying the condition of pseudo-passivity. In particular, the following theorem holds.

\begin{thm}\label{thm:implications}
Let $R = Y*$ be a convolution operator satisfying the condition of pseudo-passivity with $\vec{d} = \vec{0}$. Then, $R$ also satisfies the conditions of causality and slow growth.
\end{thm}

The proof of this theorem is split into two parts, covered by Propositions \ref{prop:implication_of_causality} and \ref{prop:implication_of_slow_growth}, respectively.

\subsection{Causality}\label{subsec:causality}

That passivity implies causality was first observed by Youla, Castriota and Carlin \cite[Thm. 1]{YoulaEtal1959}, see also \cite[Lem. 3, pg. 303]{Zemanian1965}. In order to establish whether causality is also implied by the condition of pseudo-passivity, we recall a lemma that establishes that the condition of causality needs only to be checked on the subset $\iota(\dd) \subseteq \Dom(R)$ for it to hold on all of $\Dom(R)$.

\begin{lemma}\label{lem:checking_causality}
\cite[Lem. 2, pg. 301]{Zemanian1965} Let $R$ be a convolution operator and write $\psi := \iota^{-1}(R(\iota(\varphi)))$ for $\varphi \in \dd$. If for any number $t_0 \in \R$, we have that $\psi(\xi) = 0$ for $\xi \in (-\infty,t_0)$ whenever $\varphi(\xi) = 0$ for $\xi \in (-\infty,t_0)$, then the operator $R$ satisfies the condition of causality.
\end{lemma}

We may now adopt the proof of the classical result \cite[Lem. 3, pg. 303]{Zemanian1965} in order to prove the first part of Theorem \ref{thm:implications}.

\begin{prop}\label{prop:implication_of_causality}
Let $R = Y*$ be a convolution operators satisfying the condition of pseudo-passivity with $\vec{d} = \vec{0}$. Then, $R$ also satisfies the condition of causality.
\end{prop}

\proof
Let $\varphi_1 \in \dd$ be an arbitrary test function with $\psi_1 := \iota^{-1}(R(\iota(\varphi_1)))$. Let also $a \in \R$ be a free parameter. Take now $t_0 \in \R$ arbitrary and let $\varphi \in \dd$ be a test function, such that $\varphi(\xi) = 0$ for all $\xi \in (-\infty,t_0)$. If we manage to show that the function $\psi := \iota^{-1}(R(\iota(\varphi)))$ is also identically zero in the same interval, then the result follows by Lemma \ref{lem:checking_causality}.

Define now a new test function $\varphi_2 := \varphi_1 + a\varphi$. Then, it holds for $\psi_2 := \iota^{-1}(R(\iota(\varphi_2)))$ that $\psi_2 = \psi_1 + a\psi$. Furthermore, by assumption on the operator $R$, it holds for any $x < t_0$ that
$$\Re\left[\int_{-\infty}^x\bar{\varphi_2(\tau)}\psi_2(\tau)\diff\tau\right] \geq \sum_{j = 0}^{N}c_j\int_{-\infty}^x|\varphi_2^{(j)}(\tau)|^2\diff\tau$$
for some number $N \in \N_0$ and some vector $\vec{c} \in \R^{N+1}$. Now, using the property that $\varphi(\xi) = 0$ for all $\xi \in (-\infty,t_0)$, we calculate that
\begin{multline*}
\Re\left[\int_{-\infty}^x\bar{\varphi_2(\tau)}\psi_2(\tau)\diff\tau\right] = \Re\left[\int_{-\infty}^x\bar{(\varphi_1+a\varphi)(\tau)}(\psi_1 + a\psi)(\tau)\diff\tau\right] \\
= \Re\left[\int_{-\infty}^x\bar{\varphi_1(\tau)}\psi_1(\tau)\diff\tau\right] + a\:\Re\left[\int_{-\infty}^x\bar{\varphi_1(\tau)}\psi(\tau)\diff\tau\right].
\end{multline*}

On the other hand, due to the same property of the function $\varphi$, we also have
$\varphi_2(\xi) = \varphi_1(\xi)$ for all $\xi \in (-\infty,x]$, which also transfers over to all derivatives of the functions $\varphi_1$ and $\varphi_2$, \ie $\varphi_2^{(j)}(\xi) = \varphi_1^{(j)}(\xi)$ for all $\xi \in (-\infty,x]$ and all $j \in \N$. Thus, we conclude that
$$\Re\left[\int_{-\infty}^x\bar{\varphi_1(\tau)}\psi_1(\tau)\diff\tau\right] + a\:\Re\left[\int_{-\infty}^x\bar{\varphi_1(\tau)}\psi(\tau)\diff\tau\right] \geq \sum_{j = 0}^{N}c_j\int_{-\infty}^x|\varphi_1^{(j)}(\tau)|^2\diff\tau,$$
where the inequality holds for any value of the free parameter $a$. But this now means that
$$\Re\left[\int_{-\infty}^x\bar{\varphi_1(\tau)}\psi(\tau)\diff\tau\right] = 0$$
for all $x \in (-\infty,t_0)$, as, otherwise, the parameter $a$ could be chosen large enough (in absolute value) as to break the inequality.

Observe now that the last equality may be rewritten as
\begin{multline*}
\Re\left[\int_{-\infty}^x\bar{\varphi_1(\tau)}\psi(\tau)\diff\tau\right] \\
= \int_{-\infty}^x\Re[\varphi_1](\tau)\Re[\psi](\tau)\diff\tau + \int_{-\infty}^x\Im[\varphi_1](\tau)\Im[\psi](\tau)\diff\tau = 0.
\end{multline*}
Since $\varphi_1 \in \dd$ is arbitrary, the result follows.
\endproof

\begin{remark}
If $\vec{d} \neq \vec{0}$, the above proof fails as it does not necessarily hold that
$$\sum_{j = 0}^{N}d_k\int_{-\infty}^x|\psi_2^{(k)}(\tau)|^2\diff\tau \geq \sum_{j = 0}^{N}d_k\int_{-\infty}^x|\psi_1^{(k)}(\tau)|^2\diff\tau.$$
\end{remark}

\subsection{Slow growth}\label{subsec:slow_growth}

That passivity implies slow-growth was observed, for example, by Zemanian \cite{Zemanian1963}. In fact, his result is slightly stronger as it only requires that the operator $R = Y*$ satisfies the condition of weak passivity on $\iota(\dd)$, \ie for any test function $\varphi \in \dd$ it holds that $R(\iota(\varphi)) \in \iota(\Loc)$ and
\begin{equation}\label{eq:weak_passivity}
\Re\left[\int_{-\infty}^\infty\bar{\varphi(\tau)}\psi(\tau)\diff\tau\right] \geq 0,
\end{equation}
where $\psi := \iota^{-1}(R(\iota(\varphi)))$. This turns out to be sufficient to conclude that the defining distribution $Y$ is non-negative definite, \cf \cite[Thm. 1]{Zemanian1963}, and afterwards invoking a result of Schwartz saying that any non-negative definite distribution is in $\sd'$ \cite{Schwartz1959}.

We may now prove the second part of Theorem \ref{thm:implications}.

\begin{prop}\label{prop:implication_of_slow_growth}
Let $R = Y*$ be a convolution operators satisfying the condition of pseudo-passivity with $\vec{d} = \vec{0}$. Then, $R$ also satisfies the condition of slow growth.
\end{prop}

\proof
Due to the assumption on the operator $R$, it holds, in particular, that
\begin{equation}\label{eq:weak_non_passivity_v2}
\Re\left[\int_{-\infty}^\infty\bar{\varphi(\tau)}\psi(\tau)\diff\tau\right] \geq \sum_{j = 0}^{N}c_j\int_{-\infty}^\infty|\varphi^{(j)}(\tau)|^2\diff\tau
\end{equation}
for any $\varphi \in \dd$, with $\psi := \iota^{-1}(R(\iota(\varphi)))$. Furthermore, observe that for any $j \in \N_0$, it holds that
\begin{multline*}
\int_{-\infty}^\infty|\varphi^{(j)}(\tau)|^2\diff\tau \\ 
= (-1)^j\int_{-\infty}^{\infty}\bar{\varphi(\tau)}\varphi^{(2j)}(\tau)\diff\tau = (-1)^j\int_{-\infty}^{\infty}\bar{\varphi(\tau)}\iota^{-1}(\delta_0^{(2j)} * T_\varphi)(\tau)\diff\tau.
\end{multline*}
Here, the first of the above equalities follows after integrating by parts $j$-times, while the second equality holds due to the fact that
$$\delta_0^{(m)} * T_\varphi = T_{\varphi^{(m)}}$$
for any $\varphi \in \dd$ and any $m \in \N_0$. 

We conclude now that inequality \eqref{eq:weak_non_passivity_v2} may be written as
$$\Re\left[\int_{-\infty}^\infty\bar{\varphi(\tau)}\iota^{-1}(\til{Y}*T_\varphi)(\tau)\diff\tau\right] \geq 0,$$
where the distribution $\til{Y}$ is defined as
$$\til{Y} := Y - \sum_{j = 0}^{N}(-1)^jc_j\delta_0^{(2j)}.$$
In other words, we have shown that the operator $\til{R} := \til{Y}*$ satisfies the condition of weak passivity on $\iota(\dd)$, \cf condition \eqref{eq:weak_passivity}.

By the discussion at the beginning of Section \ref{subsec:slow_growth}, it now holds that $\til{Y} \in \sd'$, yielding further that $Y \in \sd'$. This finishes the proof.
\endproof

\begin{remark}
If $\vec{d} \neq \vec{0}$, then the above method of proof fails as we are unable to write
$$\int_{-\infty}^\infty|\psi^{(j)}(\tau)|^2\diff\tau = k_j \int_{-\infty}^{\infty}\bar{\varphi(\tau)}\iota^{-1}(U * T_\varphi)(\tau)\diff\tau$$
for some $k_j \in \C$ and some distribution $U \in \sd'$.
\end{remark}

\section{Scattering pseudo-passive operators}\label{sec:scattering_operators}

The theory of (admittance) passive operators comes with a related theory of \emph{scattering passive operators}, which may be defined as follows, \cf \cite{WohlersBeltrami1965,Zemanian1968a}.

\begin{define}\label{def:scattering_passive_operator}
A convolution operator $S = Z *$ is called a \emph{scattering passive causal operator of slow growth} if it satisfies the following three conditions:
\begin{itemize}
    \item[(s)]{\emph{scattering passivity} on $\iota(\dd)$: for any number $t \in \R$ and any test function $\zeta \in \dd$ it holds that $S(\iota(\zeta)) \in \iota(\dd)$ and
\begin{equation}\label{eq:scattering_passivity}
\int_{-\infty}^t\left(|\zeta(\tau)|^2-|\eta(\tau)|^2\right)\geq 0,
\end{equation}
where $\eta := \iota^{-1}(S(\iota(\zeta)))$,}
    \item[(b)]{\emph{causality}: $\supp(Z) \subseteq [0,\infty)$,}
    \item[(c)]{\emph{slow growth}: $Z \in \sd'$.}
\end{itemize}
\end{define}

Inspired by Definition \ref{def:non_passive_operator}, we may now consider the following generalization of scattering passive operators.

\begin{remark}
Below, the Kronecker $\delta$-symbol is written as $\di_{0,j}$ instead of $\delta_{0,j}$ in order to avoid confusion with the Dirac distribution $\delta_0$.
\end{remark}

\begin{define}\label{def:scattering_non_passive_operator}
A convolution operator $S = Z *$ is called a \emph{scattering pseudo-passive causal operator of slow growth} if it satisfies the following three conditions:
\begin{itemize}
    \item[(s')]{\emph{scattering pseudo-passivity} on $\iota(\dd)$: There exist a number $N \in \N_0$ and vectors $\vec{F},\vec{G} \in \R^{N+1}$, such that for any number $t \in \R$ and any test function $\zeta \in \dd$ it holds that $S(\iota(\zeta)) \in \iota(\dd)$ and
\begin{multline}\label{eq:scattering_non_passivity}
\sum_{j = 0}^N\int_{-\infty}^t\left((\di_{0,j}-F_j)|\zeta^{(j)}(\tau)|^2-(\di_{0,j}+F_j)|\eta^{(j)}(\tau)|^2\right)\diff\tau \\
\geq 2\sum_{j = 0}^NG_j\:\Re\left[\int_{-\infty}^t\bar{\zeta^{(j)}(\tau)}\eta^{(j)}(\tau)\diff\tau\right],
\end{multline}
where $\eta := \iota^{-1}(S(\iota(\zeta)))$ and $\di_{0,j}$ denotes the Kronecker $\delta$-symbol,}
    \item[(b)]{\emph{causality}: $\supp(Z) \subseteq [0,\infty)$,}
    \item[(c)]{\emph{slow growth}: $Z \in \sd'$.}
\end{itemize}
\end{define}

\begin{example}
Let us consider the operators $S_{\pm} := \pm\delta_0 *$. If $\zeta \in \dd$, then $S_\pm * T_\zeta = \pm T_\zeta$, yielding that 
$$\eta := \iota^{-1}(S_{\pm}(\iota(\zeta))) = \pm\zeta.$$
Using this, we calculate that
$$\int_{-\infty}^t\left((1-F_0)|\zeta^{(j)}(\tau)|^2-(1+F_0)|\zeta^{(j)}(\tau)|^2\right)\diff\tau = -2F_0\int_{-\infty}^t|\zeta(\tau)|^2\diff\tau.$$

Therefore, the operator $S_+$ satisfies condition \eqref{eq:scattering_non_passivity} for $N = 0$ and any numbers $F_0$ and $G_0$ such that $F_0 \leq -G_0$, and satisfies condition \eqref{eq:scattering_non_passivity} with equality for $N = 0$ and $F_0 = - G_0$. On the other hand, the operator $S_-$ satisfies condition \eqref{eq:scattering_non_passivity} for $N = 0$ and any numbers $F_0$ and $G_0$ such that $F_0 \leq G_0$, and satisfies condition \eqref{eq:scattering_non_passivity} with equality for $N = 0$ and $F_0 = G_0$.\hfill$\lozenge$
\end{example}

\begin{example}
Assume that an operator $S$ satisfies condition \eqref{eq:scattering_non_passivity} with equality for some number $N \in \N_0$ and vectors $\vec{F},\vec{G} \in \R^{N+1}$. Then, the operator $-S$ satisfies condition \eqref{eq:scattering_non_passivity} with equality for the same number $N$ and vectors $\vec{F},-\vec{G} \in \R^{N+1}$.\hfill$\lozenge$
\end{example}

\subsection{Transition between admittance and scattering operators}\label{subsec:transition}

The background to the transition between admittance and scattering passive operators is the algebraic equivalence between inequalities \eqref{eq:passivity} and \eqref{eq:scattering_passivity}. Indeed, one sees this by writing $\varphi = \zeta + \eta$ and $\psi = \zeta - \eta$, for $\zeta,\eta \in \dd$. Furthermore, by using the same transformation, one can also establish an algebraic equivalence between inequalities \eqref{eq:admittance_non_passivity} and \eqref{eq:scattering_non_passivity}, where the relations between the vectors $\vec{c},\vec{d}$ and the vectors $\vec{F}, \vec{G}$ become $F_j = c_j + d_j$ and $G_j = c_j - d_j$. For the transformation $\varphi = \zeta + \eta$ and $\psi = \zeta - \eta$ to make sense in terms of convolution operators, we must consider the convolution algebra of right-sided distributions.

The space $\dd'_\mathrm{r}$, \cf Section \ref{sec:background}, equipped with the operations $+$ and $*$, \ie addition and convolution of distributions, becomes an algebra over $\C$. The unity of the of convolution is the Dirac distribution, \ie
$$U * \delta_0 = \delta_0 * U = U$$
for any $U \in \dd'_\mathrm{r}$ (this holds even for $U \in \dd'$). Therefore, the convolution inverse of a distribution $U \in \dd'_\mathrm{r}$, if it exists, is a distribution $U^{*-1} \in \dd'_\mathrm{r}$, such that
$$U * U^{*-1} = U^{*-1} * U = \delta_0.$$

In our case, we must, for an admittance passive operator $R = Y*$, be able to solve the equation 
$$\zeta - \eta = \iota^{-1}(Y * T_{\zeta + \eta})$$
for the function $\eta$. Formally, the solution is
$$\eta = \iota^{-1}((\delta_0 + Y)^{*-1} * (\delta_0 - Y) * T_{\zeta}),$$
and for this expression to be well-defined, we must be able to define the convolution inverse of the distribution $\delta_0 + Y$.

If the operator $R = Y*$ is a real admittance passive operator, then we may define the distribution $(\delta_0 + Y)^{*-1}$ in the following way \cite[pg. 418]{Zemanian1968a}. Due to the assumption on $R$, the function $W := \mathcal{L}(Y)$ is a positive-real function, and the same holds for the function $s \mapsto (1 + W(s))^{-1}$. As the inverse Laplace transform of this function exists by Theorem \ref{thm:passive_operators}, we may define
$$(\delta_0 + Y)^{*-1} := \mathcal{L}^{-1}(s \mapsto (1 + W(s))^{-1}).$$
Furthermore, by Theorem \ref{thm:passive_operators}, the operator
\begin{equation}\label{eq:auxiliary_operator}
\til{R} = (\delta_0 + Y)^{*-1}*
\end{equation}
will also be a real admittance passive operator and, thus, it holds that $\supp((\delta_0 + Y)^{*-1}) \subseteq [0,\infty)$ and $(\delta_0 + Y)^{*-1} \in \sd'$. Thus, the convolution operator $S = Z*$ with
$$Z := (\delta_0 + Y)^{*-1} * (\delta_0 - Y)$$
is well defined and is a scattering passive operator, as the properties of $Y$ and $(\delta_0 + Y)^{*-1}$ guarantee that $\supp(Z) \subseteq [0,\infty)$ and $Z \in \sd'$, while the condition of scattering passivity is fulfilled automatically due to the algebraic equivalence of inequalities \eqref{eq:passivity} and \eqref{eq:scattering_passivity} as discussed previously.

If $R = Y*$ is an admittance passive operator instead, \ie the operator is no longer assumed to be real, the distribution $(\delta_0 + Y)^{*-1}$ may still be defined via the inverse Laplace transform as before. However, the operator $\til{R}$, defined as before, will not necessarily be a passive operator and the operator $S = Z*$, with the distribution $Z$ defined as before, will automatically only satisfy the condition of scattering passivity. The condition of scattering passivity does, however, independently imply that $\supp(Z) \subseteq [0,\infty)$ \cite[pg. 425]{Zemanian1968a}.

If $R = Y*$ is an admittance pseudo-passive operator, then we are, in general, unable to define the distribution $(\delta_0 + Y)^{*-1}$ via the inverse Laplace transform as was done in the previous cases. However, if it exists, we may define the operator $S = Z*$ with the distribution $Z$ defined as before, which will automatically satisfy the condition of scattering pseudo-passivity. We will discuss later, in Sections \ref{subsec:scat_causality} and \ref{subsec:scat_slow_growth}, whether the condition of scattering pseudo-passivity implies the conditions of causality and/or slow-growth for at least some special cases. Separately, the precise relationship between the existence of a Laplace transform and the condition of causality, without assuming any form of (pseudo)-passivity, has been considered in \cite{Zemanian1968b}. 

\subsection{The Laplace transform}
We may determine the restrictions on the range of the Laplace transform of the defining distribution of a scattering pseudo-passive operator of slow growth $S$ via an analogous procedure to the one presented in Section \ref{sec:Laplace}. Therefore, let $t \in \R$ be arbitrary, take any $s \in \C_+ := \{z \in \C~|~\Re[z] > 0\}$ and choose a function $\zeta \in \sd$ with the property that $\zeta(\xi) = \E^{s\xi}$ for all $\xi \in (-\infty,a)$, where $a > t$. Then, we have that $S(\iota(\varphi)) \in \iota(\sd)$ and it holds for $\eta := \iota^{-1}(S(\iota(\zeta)))$ that
$$\eta(\xi) = \langle Z, t' \mapsto \zeta(\xi - t')\rangle.$$
For the derivatives and anti-derivatives of the functions $\zeta$ and $\eta$, it holds that
$$\zeta^{(j)}(\xi) = s^j\E^{s\xi} \quad\text{and}\quad \eta^{(j)}(\xi) = s^j\E^{s\xi}W(s)$$
for any $j \in \N_0$ and any $\xi \in (-\infty,t]$, where $W := \mathcal{L}(Z)$ denotes the Laplace transform of the defining distribution $Z$ of the operator $S$.

We calculate now that
\begin{multline*}
\sum_{j = 0}^N\int_{-\infty}^t\left((\delta_{0,j}-F_j)|\zeta^{(j)}(\tau)|^2-(\delta_{0,j}+F_j)|\eta^{(j)}(\tau)|^2\right)\diff\tau \\
= \sum_{j = 0}^N|s|^{2j}\left((\delta_{0,j}-F_j)-(\delta_{0,j}+F_j)|W(s)|^2\right)\int_{-\infty}^t\E^{2\Re[s]\tau}\diff\tau
\end{multline*}
and
$$\sum_{j = 0}^NG_j\:\Re\left[\int_{-\infty}^t\bar{\zeta^{(j)}(\tau)}\eta^{(j)}(\tau)\diff\tau\right] = \sum_{j = 0}^NG_j\:|s|^{2j}\Re[W(s)]\int_{-\infty}^t\E^{2\Re[s]\tau}\diff\tau,$$
yielding that
\begin{equation}\label{eq:range_of_scattering_W}
\sum_{j = 0}^N|s|^{2j}\left((\delta_{0,j}-F_j)-(\delta_{0,j}+F_j)|W(s)|^2\right) \geq 2\sum_{j = 0}^NG_j\:|s|^{2j}\Re[W(s)].
\end{equation}
As before, if the operator $S$ was real, this translates to the Laplace transform $W$ as the property that $W(s) \in \R$ if $s \in (0, \infty)$. We summarize this result in the following proposition.

\begin{prop}\label{prop:range_of_scattering_W}
Let $S = Z*$ be a scattering pseudo-passive causal operator of slow-growth satisfying condition \eqref{eq:scattering_non_passivity} for some $N \in \N_0$ and $\vec{F},\vec{G} \in \R^{N+1}$. Then, the Laplace transform $W := \mathcal{L}(Z)$ of its defining distribution exists and is a \holo function on $\C_+$ satisfying, for each $s \in \C_+$, the inequality \eqref{eq:range_of_scattering_W}. Furthermore, if the operator $S$ is assumed to be real, it holds that $W(s) \in \R$ if $s \in (0, \infty)$.
\end{prop}

\begin{example}
As in Example \ref{ex:cd_non_passive}, we investigate in detail the geometric shape of the range of the function $W$ when $N = 0$. Using inequality \eqref{eq:range_of_scattering_W}, \ie
$$(1-F)-(1+F)|W(s)|^2 \geq 2G\:\Re[W(s)],$$
we have that $\Ran(W) \subseteq B$, where the set $B$ is defined as
\begin{multline*}
B := \{\sigma \in \C~|~(1-F)-(1+F)|\sigma|^2 \geq 2G\:\Re[\sigma]\} \\
= \{(x,y) \in \R^2~|~(1-F)-(1+F)(x^2+y^2) \geq 2G\:x\}.
\end{multline*}
The geometric picture of this set $B$ is the following.

If $F = -1$, there are no square terms in the definition of the set $B$, meaning that we get a half-plane. Else, if $F \neq -1$, we may divide by the number $1+F$ to get
$$\frac{1-F}{1+F} \geq x^2+ y^2 + \frac{2G}{1+F}x \quad\text{or}\quad \frac{1-F^2+G^2}{(1+F)^2} \geq \left(x-\frac{G}{1+F}\right)^2 + y^2$$
if $F > -1$ and
$$\frac{1-F}{1+F} \leq x^2+ y^2 + \frac{2G}{1+F}x \quad\text{or}\quad \frac{1-F^2+G^2}{(1+F)^2} \leq \left(x-\frac{G}{1+F}\right)^2 + y^2$$
if $F < - 1$. Thus, the shape of the set $B$ is characterized in terms of the parameters $F$ and $G$ in the following way.
\begin{itemize}
\item[(i)]{We have $B = \C$ if $F < -1$ and $1- F^2 + G^2 \leq 0$.}
\item[(ii)]{We have $B$ being the area outside and including the circle
\begin{equation*}\label{eq:scattering_circle}
\frac{1-F^2+G^2}{(1+F)^2} = \left(x-\frac{G}{1+F}\right)^2 + y^2
\end{equation*}
if $F < -1$ and $1 - F^2 + G^2 > 0$.}
\item[(iii)]{We again have $B = \C$ if $F = -1$ and $G = 0$.}
\item[(iv)]{We have $b$ equal to the half-plane $\{(x,y) \in \R~|~2 \geq G\:x\}$ if $F = -1$ and $G \neq 0$.}
\item[(v)]{We have $B$ being the area inside and including the circle \eqref{eq:scattering_circle} if $F > -1$ and $1 - F^2 + G^2 > 0$.}
\item[(vi)]{We have $B = \{\tfrac{G}{1+F}\} \subseteq \C$ if $F > -1$ and $1 - F^2 + G^2 = 0$.}
\item[(vii)]{We have $B$ being empty if $F > -1$ and $1 - F^2 + G^2 < 0$. \hfill$\lozenge$}
\end{itemize}
\end{example}

\subsection{Causality}\label{subsec:scat_causality}

For scattering pseudo-passive causal operators of slow growth, we may establish two instances where the condition of scattering pseudo-passivity implies the condition of causality.

\begin{prop}\label{prop:implication_of_scat_causality_v1}
Let $S = Z*$ be a convolution operator satisfying the condition of scattering pseudo-passivity \eqref{eq:scattering_non_passivity} with $F_0 \geq -1$, $F_j \geq 0$ for $j > 0$ and $\vec{F} = \vec{G}$. Then, $S$ also satisfies the condition of causality.
\end{prop}

\proof
Take $t_0 \in \R$ arbitrary and let $\zeta \in \dd$ be a test function, such that $\zeta(\xi) = 0$ for all $\xi \in (-\infty,t_0)$. By Lemma \ref{lem:checking_causality}, it suffices to show that the function $\eta := \iota^{-1}(S(\iota(\zeta)))$ is also identically zero in the same interval.

Taking into account the assumption that $\vec{F} = \vec{G}$, we may rewrite condition \eqref{eq:scattering_non_passivity} as
$$\int_{-\infty}^x|\zeta(\tau)|^2\diff\tau \geq \int_{-\infty}^x|\eta(\tau)|^2\diff\tau + \sum_{j = 0}^NF_j\int_{-\infty}^x|(\zeta^{(j)}-\eta^{(j)})(\tau)|^2\diff\tau,$$
where $x \in (-\infty,t_0)$. Using also the assumption that $\zeta(\xi) = 0$ for all $\xi \in (-\infty,t_0)$, the previous inequality simplifies to
$$0 \geq \sum_{j = 0}^N(\delta_{0,j} + F_j)\int_{-\infty}^x|\eta^{(j)}(\tau)|^2\diff\tau.$$
However, by the assumptions on the entries of the vector $\vec{F}$, the right-hand side of the above inequality is non-negative, finishing the proof.
\endproof

\begin{prop}\label{prop:implication_of_scat_causality_v2}
Let $S = Z*$ be a convolution operator satisfying the condition of scattering pseudo-passivity \eqref{eq:scattering_non_passivity} with $\vec{F} = \vec{0}$. Then, $S$ also satisfies the condition of causality.
\end{prop}

\proof
By the assumption on the operator $S$, we may, for functions $\zeta$ and $\eta$ as in the proof of Proposition \ref{prop:implication_of_scat_causality_v1}, rewrite condition \eqref{eq:scattering_non_passivity} as
$$\int_{-\infty}^x|\zeta(\tau)|^2\diff\tau \geq \int_{-\infty}^x|\eta(\tau)|^2\diff\tau + 2\sum_{j = 0}^NG_j\:\Re\left[\int_{-\infty}^x\bar{\zeta^{(j)}}(\tau)\eta^{(j)}(\tau)\diff\tau\right].$$
Taking into account the assumption that $\zeta(\xi) = 0$ for all $\xi \in (-\infty,t_0)$ yields
$$0 \geq \int_{-\infty}^x|\eta(\tau)|^2\diff\tau,$$
which, by Lemma \ref{lem:checking_causality}, finishes the proof.
\endproof

\subsection{Slow growth} \label{subsec:scat_slow_growth}

Whether the condition of scattering pseudo-passivity implies the condition of slow growth remains inconclusive. The proof of Proposition \ref{prop:implication_of_slow_growth} relied on the possibility to construct, out of the defining distribution $Y$ of the operator $R$, a distribution $\til{Y}$ such that the operator $\til{R} = \til{Y}*$ satisfied the condition of weak passivity and, afterwards, relying on Zemanian's result \cite[Thm. 1]{Zemanian1963}. Already for a convolution operator satisfying the condition of scattering passivity, such a method of proof fails due to the absence of an analogue of Zemanian's result for the scattering case.

\section{Conclusion}\label{sec:summary}

We infer from Theorems \ref{thm:range_of_W} and \ref{thm:implications} that pseudo-passive causal operators of slow growth satisfying condition \eqref{eq:admittance_non_passivity} with $\vec{d} = \vec{0}$ exhibit many of the properties that hols for the classic case of passive operators. Namely, the condition of pseudo-passivity with $\vec{d} = \vec{0}$ implies the conditions of causality and slow growth, and a restriction on the range of the Laplace transform of the defining distribution of the operator may be obtained. For operators which satisfy condition \eqref{eq:admittance_non_passivity} with $\vec{d} \neq \vec{0}$, we are still able to obtain a restriction on the range of the Laplace transform, but we are, using the methods of proof presented here, unable to determine whether the conditions of causality and slow growth are still automatically implied.

For scattering pseudo-passive causal operators of slow growth, restrictions on the range of the Laplace transform of their defining distributions are also obtained, \cf Proposition \ref{prop:range_of_scattering_W}, but the relations between the conditions scattering pseudo-passivity, causality and slow-growth are not strong as in the previous case, \cf Propositions \ref{prop:implication_of_scat_causality_v1} and \ref{prop:implication_of_scat_causality_v2}. 

\section*{Acknowledgements}

The author would like to thank Mats Gustafsson for initiating this work by bringing up the topic, in particular for conjecturing the restrictions on the range of the Laplace transform of convolution operators satisfying condition \eqref{eq:admittance_non_passivity} for $N = 0$ and $d_0 = 0$ and of convolution operators satisfying condition \eqref{eq:scattering_non_passivity} for $N = 0$ and $G_0 = 0$.

Furthermore, the author would like to thank Odysseas Bakas, Dale Frymark and Annemarie Luger for many helpful comments and careful reading of the manuscript.

\bibliographystyle{amsplain}
\bibliography{total,total2}

\end{document}